\documentclass[12pt,reqno]{article}
\usepackage{amsfonts}
\usepackage{graphicx}
\usepackage{epstopdf}
\usepackage{amsmath}
\usepackage{amsthm}
\usepackage{latexsym,bm}
\usepackage{amssymb}
\usepackage{indentfirst}
\newtheorem{thm}{Theorem}[section]

\newtheorem{lem}[thm]{Lemma}

\numberwithin{equation}{section}

\begin{document}

\title{{\bf
The index of isolated umbilics on surfaces of non-positive
curvature }}
\author{\bf F. Fontenele\thanks{Work partially supported by CNPq
(Brazil)}\, and F. Xavier\thanks{Work partially supported by CAPES
and IMPA (Brazil)}\\ \\\small{\it To Professor Jorge Sotomayor, on
the occasion of his 70th birthday}}

\date{}
\maketitle

\begin{quote}

\noindent \small {\bf Abstract}. It is shown   that if a $C^2$
surface $M\subset \mathbb R^3$ has negative curvature on
the complement of a point $q\in M$, then the $\mathbb Z/2$-valued Poincar\'e-Hopf  index at  $q$ of either
distribution of principal directions  on $M-\{q\}$ is
non-positive.  Conversely, any non-positive half-integer arises in this fashion. The proof of the index estimate is based  on  geometric-topological arguments, an index theorem for
symmetric  tensors on  Riemannian surfaces,  and  some aspects of the  classical Poincar\'e-Bendixson theory.

\end{quote}

\vskip20pt

\section{Introduction.}

\vskip10pt

The   distributions of principal directions   on a  surface in
$\mathbb R^3$,  defined on the complement of the umbilical set (i.e., the locus where the principal curvatures coincide),
have been the object of intense scrutiny since the early days of
differential geometry.  For  both technical and geometric reasons,
most of these investigations were  conducted  under the hypothesis
that the surface   is real analytic,  or  at least  of class
$C^3$,  although one needs only $C^2$ regularity in order for the
fields of principal directions to be continuous.

The aim of this work is to  establish
 an   estimate for the local index of these
fields,
under a natural  curvature restriction,   but  with optimal regularity:

\begin{thm} \label{teorema 1}
Let $M\subset\mathbb R^3$ be a $C^2$ surface, $q\in M$. Assume
that the Gaussian curvature of $M$ is negative at every point of $M$ other than $q$.
Then, the $\mathbb Z/2$-valued Poincar\'e-Hopf index at
$q$ of either distribution of principal directions on $M-\{q\}$ is
non-positive.
\end{thm}
We observe that the theorem is sharp. Indeed, if the Gaussian curvature remains negative at $q$,  then the distributions of
principal directions extend continuously to $M$, and so the index
is zero. On the other hand,  given  any     negative  number
$j\in  \mathbb Z/2$  one can construct a surface as in the
statement of the theorem, even a minimal one (i.e., with vanishing mean curvature), that has an isolated
umbilic of index  $j$ (\cite{SX3}).

When the surface in question is minimal, the conclusion in Theorem \ref{teorema 1} can be verified
using the holomorphic
data  in the
Weierstrass representation of the surface.

The main point of the present work is that,  surprisingly,  Theorem \ref{teorema 1}
 applies  to surfaces that are merely $C^2$,  and not only to those    $C^{\omega}$ surfaces that are minimal. Since  complex analysis is no longer available in this more general setting,   new tools have to be introduced in order  to estimate  the index. Loosely speaking, the replacement for complex analysis is, when properly augmented,  the classical qualitative theory of planar dynamical systems.

Umbilic points are notoriously elusive geometric objects. For instance,  on a surface of positive curvature the index of an isolated umbilic need not be positive. In other words,   the ``dual" statement of Theorem
\ref{teorema 1} does not hold.  Indeed,  inverting
surfaces  satisfying the hypotheses of Theorem \ref{teorema 1} on   suitable spheres -- a process that does not change
the index of the umbilic --, one can produce surfaces of \it positive \rm
Gaussian curvature  exhibiting an umbilic  whose index is any
prescribed  \it negative \rm   half-integer.

By analogy with the above mentioned sharpness of Theorem
\ref{teorema 1},  one might naively expect that every \it positive \rm  half-integer could be realized as the index of an isolated umbilic on a surface of \it positive \rm  curvature. In stark opposition to  this expectation,  it is actually predicted that
on any  sufficiently regular surface, without any curvature
restrictions whatsoever,  the index of an isolated umbilic should be  at most \it one.  \rm
This is the  well-known   local Carath\'eodory conjecture, also
known as the Loewner conjecture. We refer the reader to
\cite{Bo}, \cite{Br}, \cite{FB}-\cite{H2},  \cite{IV}-\cite{X} for
a sample of the many  works that have appeared in print, old and new,  on this very challenging problem,  as well as  on
various aspects of the global study of principal foliations.

To put our results in the  context of  the Carath\'eodory conjecture, we re-iterate that Theorem \ref{teorema
1}  is sharp and verifies the $C^2$ version of the local Carath\'eodory conjecture in a
 geometrically important special case, \it but with  a  stronger
conclusion. \rm
Thus,   Theorem \ref{teorema
1}  represents a contribution to the interface between  classical differential geometry and classical dynamical systems  that   stands on its own, since it  cannot be subsumed by the  resolution of the Carath\'eodory conjecture.  On the other hand, we hasten to add that  there is no expectation that the present method can be used to tackle the said  conjecture,  given our strong reliance on negative curvature.

Over the years, the task of
estimating the index of isolated umbilics has proven to be an
arduous   one,   often  involving lengthy  and intricate arguments (e.g., \cite{IV}).
Against this backdrop,
it is pleasing that the  proof of Theorem
\ref{teorema 1}, albeit  delicate in its own right,   is  rather  conceptual.  The arguments are
based on an index  theorem  for  abstract symmetric tensors on  Riemannian
surfaces,  elements of  the classical qualitative
theory of two-dimensional dynamical systems, and a modicum of topology and classical differential geometry.
Although the main result is new, the subject matter lends itself  to a more expository style and, accordingly,   full details are provided.

\section {An index theorem for abstract symmetric  tensors.}

\vskip10pt  A continuous symmetric tensor field $A$ of type ( 1,1),
defined on  a Riemannian surface $M$, is said to have  an isolated
$A$-singularity at $q\in M$ if  there exists a neighborhood  $V$
of $q$ such that the  eigenvalues of $A(p)$ are unequal   for any
$p\in V-\{q\}$. This condition is equivalent to the requirement
that the traceless part  $A(p)-\frac{1}{2}(\textnormal{tr} A(p))I$
of $A(p)$ be non-zero for all $p\neq q$.  If  the eigenvalues of
$A(q)$ are  distinct, then $q$ is automatically an isolated
$A$-singularity (one then thinks of $q$ as being a ``removable"
$A$-singularity).

Under the above conditions, there are two continuous line fields
on $M-\{q\}$, not necessarily orientable,  which correspond  to
the diagonalizing directions of $A$.   Given  a continuous field
of directions $\xi$ on $V-\{q\}$,  denote  by $A\xi$ the field of
directions obtained by  applying $A(p)$ to any vector generating
the one-dimensional subspace $\xi(p)$ of $T_pM$. We write  $j(A)$
for the index at $q$ of either of the two  fields of diagonalizing
directions of $A$. Similarly,  $j(\eta)$   stands for  the index
at $q$ of a line field $\eta$ with an isolated singularity at $q$.
  \vskip10pt
\begin{thm}\label{FormInd} Let $A$ be a continuous symmetric tensor field of type
(1,1) on a Riemannian 2-manifold $M$, and $q\in M$ an  isolated $A$-singularity.  Then, for
every continuous field of directions $\xi$ on a punctured
neighborhood of $q$,  one has
\begin{eqnarray}\label{Field1}
2\,j(A)=j\big((A-\frac{1}{2}(\textnormal{tr}
A)I)\,\xi\big)+j(\xi).
\end{eqnarray}
\end{thm}
\vskip10pt
\noindent
\bf Remarks. \rm One should think of $\xi$ as being a ``test" line field. For instance, if $\xi$ is chosen to be one of
the continuous fields of eigendirections of $A$, each term on the
right hand side of (\ref{Field1}) equals $j(A)$. As another
illustration to see that $2$ is the correct  factor in the left
hand side, let $A$ be a continuous symmetric tensor field on a
compact orientable surface $M$ of non-zero genus,   with the
property that the set  $F_1$ where   $A$ is a multiple of the
identity is finite.  Let $\xi$ be a continuous line  field on $M$
with  a finite set $F_2$ of  singularities. Applying
(\ref{Field1}) around each point in $F_1\cup F_2$,  and summing,
one sees from the Poincar\'e-Hopf theorem that both sides of
(\ref{Field1}) equal $2\chi (M)$.

\vskip10pt

\noindent{\bf Proof.} Formula (\ref{Field1}) was first   established in \cite{SX3},  under the more restrictive assumptions that $A$ is a smooth tensor and  $\xi$  is a smooth vector field; here, $A$ and $\xi$ are  assumed to be only continuous, and $\xi$ is allowed to be an unorientable  line field. This extra generality requires a new line of argument.

 Let $\lambda\geq\mu$ be the eigenvalues of
$A$ and $U$ a neighborhood of $q$  such that $\lambda(p)>\mu(p)$  whenever  $p\in U-\{q\}$. The continuous distributions $D_{\lambda}$ and $D_{\mu}$
of eigenspaces determined by $A$ on $U-\{q\}$ have the same index
at $q$, and this common value is, by definition, the index $j(A)$
of $A$ at $q$. Let $B=A-\frac{1}{2}(\text{tr} A)I$, so that
$\text{tr}B=0$ everywhere, and let
$$\left[\begin{matrix}a&b\\b&-a\end{matrix}\right]$$ be the matrix
representation of $B$ with respect to an orthonormal frame
$\{e_1,e_2\}$ in (a possibly smaller neighborhod) $U$. Observe that the continuous vector field
$X:=ae_1+be_2$ has no zeros on $U-\{q\}$.

Let $C\subset U$ be a Jordan curve around $q$ and $\gamma:[0,1]\to
U$ a positive parametrization of $C$ (relative to the orientation
determined by $\{e_1,e_2\}$ on $U$). Let
$V(t)=\text{cos}\,\theta(t)\,e_1+\text{sin}\,\theta(t)\,e_2$ be a
continuous unit vector field along $\gamma$ such that $V(t)$
generates $D_{\lambda}(\gamma(t))$, for every $t\in [0,1]$. By
definition,
\begin{eqnarray}\label{Field2}
j(D_{\lambda})=\frac{\theta(1)-\theta(0)}{2\pi}\,\cdot
\end{eqnarray}
(Notice that since $V(1)=\pm V(0)$, $\theta(1)$ differs from
$\theta(0)$ by a multiple of $\pi$, and thus $j(D_{\lambda})\in
\frac{1}{2}\mathbb Z$). Consider the continuous vector field $W$
along $\gamma$ defined by
\begin{eqnarray}\label{Field3}
W(t)=\text{cos}(2\,\theta(t))\,e_1+\text{sin}(2\,\theta(t))\,e_2.
\end{eqnarray}
We claim that $W(t)$ is orthogonal to $(-be_1+ae_2)(\gamma(t))$,
for all $t\in [0,1]$. In fact, since $V(t)$ is an eigenvector of
$A(\gamma(t))$ (and hence  of $B(\gamma(t))$) and
$\text{sin}\,\theta(t)\,e_1-\text{cos}\,\theta(t)\,e_2$ is
orthogonal to $V(t)$, one has
\begin{eqnarray}\label{Field3a}
0&=&\big\langle
B(V(t)),\text{sin}\,\theta(t)\,e_1-\text{cos}\,\theta(t)\,e_2\big\rangle\nonumber\\&=&\big\langle
\text{cos}\,\theta(t)(ae_1+be_2)+\text{sin}\,\theta(t)(be_1-ae_2),\text{sin}\,\theta(t)\,e_1
-\text{cos}\,\theta(t)\,e_2\big\rangle\nonumber\\&=&-b\big(\text{cos}^2\theta(t)-\text{sin}^2\theta(t)\big)
+2a\,\text{sin}\,\theta(t)\,\text{cos}\,\theta(t)\nonumber\\&=&-b\,\text{cos}(2\,\theta(t))+a\,\text{sin}(2\,\theta(t))
\nonumber\\&=&\big\langle
\text{cos}(2\,\theta(t))\,e_1+\,\text{sin}(2\,\theta(t))\,e_2,-be_1+ae_2
\big\rangle\nonumber\\&=&\big\langle W(t),-be_1+ae_2 \big\rangle,
\end{eqnarray}
which proves the claim.

Since $X$ is orthogonal to $-be_1+ae_2$, it follows from the claim
above that
\begin{eqnarray}\label{Field4}
\frac{X}{|X|}\circ\gamma=\pm W.
\end{eqnarray}

Let $\xi$ be a continuous field of directions on $U-\{q\}$. If
$Z(t)=\text{cos}\,\varphi(t)\,e_1+\text{sin}\,\varphi(t)\,e_2$ is
a continuous vector field along $\gamma$ such that, for all $t\in
[0,1]$, $Z(t)$ generates $\xi(\gamma(t))$, then
\begin{eqnarray}\label{Field5a}
j(\xi)=\frac{\varphi(1)-\varphi(0)}{2\pi}\,\cdot
\end{eqnarray}
From (\ref{Field3}) and (\ref{Field4}), we obtain
\begin{eqnarray}\label{Field7}
B(Z(t)))&=&[a\text{cos}\varphi+b\text{sin}\varphi]e_1+[b\text{cos}\varphi-a\text{sin}\varphi]e_2\nonumber\\&=&
\pm|X|\big\{[\text{cos}(2\theta)\,\text{cos}\varphi+\text{sin}(2\theta)\,\text{sin}\varphi]e_1
+[\text{sin}(2\theta)\,\text{cos}\varphi-\text{cos}(2\theta)\,\text{sin}\varphi]e_2\big\}
\nonumber\\&=&\pm|X|\big\{\text{cos}\big(2\,\theta-\varphi\big)e_1+\text{sin}\big(2\,\theta-\varphi\big)e_2\big\}.
\end{eqnarray}
Since $Z(t)$ generates $\xi(\gamma(t))$, the equality above shows
that, for all $t\in [0,1]$,
\begin{eqnarray}\label{Field7a}
\text{cos}\big(2\,\theta(t)-\varphi(t)\big)e_1+\text{sin}\big(2\,\theta(t)-\varphi(t)\big)e_2
\end{eqnarray}
generates $(B\xi)(\gamma(t))$. The lemma now follows from  (\ref{Field2}) and
(\ref{Field5a}):
\begin{eqnarray*}\label{Field8}
j(B\xi)&=&\frac{[2\,\theta(1)-\varphi(1)]-[2\,\theta(0)-\varphi(0)]}{2\pi}\nonumber\\
&=&\frac{2\,\theta(1)-2\,\theta(0)}{2\pi}-\frac{\varphi(1)-\varphi(0)}{2\pi}\nonumber\\
&=&2j(D_{\lambda})-j(\xi)\nonumber\\&=&2j(A)-j(\xi).
\end{eqnarray*}

\vskip10pt

\section{Gradients and degenerate local homeomorphisms.}

\vskip10pt

\noindent  The lemma below is well known  for the usual gradient of a planar function. Here, we work in  the context of arbitrary Riemannian  surfaces.
\vskip10pt
\begin{lem}\label{IndGrad}
Let $f$ be a $C^1$ function defined on an open set $U$ of a  $C^2$ Riemannian  surface $M$, and $q\in U$ an isolated   critical point  of  $ f$.
Then, the Poincar\'e-Hopf index  of $\nabla f$ at $q$ is at most one. Furthermore,
the  index  of $\nabla f$ at $q$ is one  if and only if $f$ has a
strict local maximum, or minimum,  at $q$.
\end{lem}
\vskip10pt
\noindent{\bf Proof.} Taking  $U$ to  be a coordinate neighborhood,  $U=\varphi (W)$, $W\subset \mathbb R^2$,  we may consider a  continuous  tensor $\tilde J$ on $U$ corresponding to  rotation by $\pi/2$ in the tangent spaces of $M$.
Since $\tilde J$ can be continuously deformed into the identity through pointwise  invertible  tensors, the index at $q$ of  $\nabla f$ satisfies
$j(\nabla f)=j(\tilde J\nabla f)$. Notice that   $\tilde J\nabla f$  is tangent to  the level curves of $f$.  Likewise,  $\varphi^{\ast}(\tilde J\nabla f)$  is tangent  to the level curves of $f\circ \varphi$ on $W$. Hence
\begin{eqnarray*} j(\nabla f)=j(\tilde J\nabla f)= j(\varphi^{\ast} (\tilde J\nabla f))=j(J \nabla_0(f\circ \varphi)),
\end{eqnarray*}
where the last two indices are computed at $\varphi^{-1}(q)$, $J$ stands for the usual complex structure in $\mathbb R^2$,  and $\nabla_0$ is the Euclidean  gradient.
It  follows from  the
Poincar\'e-Bendixson theory (see, e.g.,  the exercise on p. 173 of \cite{H}) that
\begin{eqnarray*}j(J \nabla_0(f\circ \varphi)) \leq 1, \end{eqnarray*}
 with equality holding if and only if  $\varphi^{-1}(q)$  is a point of local maximum,  or minimum,  of
 $f\circ \varphi$.  In particular, $j(\nabla f)\leq 1$, and  equality holds  if and only if $q$ is an extremum of $f$.
\;\;\; q.e.d.
\vskip15pt

Under the hypotheses of Theorem \ref{teorema 1}, if $q$ is an umbilic
point then the Gaussian curvature necessarily has to vanish at $q$, and so the
Gauss map is not  a local  diffeomorphism. However, one can still
prove that the Gauss map is open. More generally, using arguments
from algebraic topology, it is possible to argue that   a
continuous map must be open if it  is a local homeomorphism on the
complement of a sufficiently ``thin" subset of  its  domain (\cite
{Ch}, p.354).  Fortunately, in the special case  that concerns us,
an  elementary proof is available:

\vskip10pt

\begin{lem}\label{OpenMap} Let $U\subset \mathbb R^n$ be open,
$n\geq 2$, $q\in U$, $F:U\to \mathbb R^n$ continuous. If the
restriction of $F$ to $U-\{q\}$ is a local homeomorphism, then $F$
is an open map. \end{lem}

\vskip5pt

\noindent{\bf Proof.} (The simple example $f(x) =x^2$   shows  the
need to have $n\geq 2$.) Assume, by contradiction, that   $F$ is
not an open map. Since the restriction of $F$ to $U-\{q\}$ is a
local homeomorphism, there exists an open set $V\subset U$, with
$q\in V$, such that $F(q)\in\partial F(V)$. Let $B$ be an open
ball centered at $q$ such that $\overline B\subset V$. We are going to need claims i) and ii) below:

\vskip5pt

\noindent i) For every $y\in\partial F(B),\;y\neq F(q)$, there
exists $x\in\partial B$ satisfying $F(x)=y$.

\vskip5pt

\noindent Indeed, let $(x_k)$ be a sequence  in $B$ with
$F(x_k)\to y$. Passing to a subsequence, we can suppose $x_k\to
x\in\overline B$. Hence $F(x_k)\to F(x)$, and so $F(x)=y$. Since
$y\neq F(q)$, and the image of every point in $B-\{q\}$ belongs to
the interior of $F(B)$, one has $x\in\partial B$.

\vskip5pt

\noindent ii) Every ball $D$ centered at $F(q)$ contains a point
$y\in\partial F(B)$ distinct from $F(q)$.

\vskip5pt

\noindent Observe that $F(q)\in \partial F(B)$, otherwise $F(q)\in
\text {int} F(B) \subset \text {int} F(V)$. Since $F(q)\in D\cap
\partial F(B)$, one cannot have $D\subset F(B)$. Let then $z\in
D-F(B)$. Using the continuity of $F$ and the fact that the
restriction of $F$ to $B-\{q\}$ is a local homeomorphism, we see
that $D\cap \text {int} F(B) \neq \emptyset.$ Since $n\geq 2$, one
can choose $w\in D\cap \text {int} F(B)$ such that $F(q)$ lies
outside the segment $\overline{wz}$ joining $w$ to $z$. Since
$\overline{wz}$ joins a point in the interior of $F(B)$ to a point
in the complement of $F(B)$, it must contain a point $y\in
\partial F(B)$, which is necessarily distinct from $F(q)$. Since
$y\in\overline{wz}\subset D$,  ii) follows.

Applying ii) to a sequence of balls $D=D_k$ centered at $F(q)$,
with radii tending to zero, one sees that there exists a sequence
$(y_k)$ in $\partial F(B)$, with $y_k\neq F(q)$ for all $k$, such
that $y_k\to F(q)$. By i), $y_k=F(x_k)$ for some sequence $(x_k)$
in $\partial B$. Passing to a subsequence, we can assume that
$x_k\to x\in\partial B$. By continuity, $F(x_k)\to F(x)$, and so
$F(x)=F(q)$. We have then found a point $x\in
\partial B \subset V-\{q\}$ whose image by $F$ belongs to
$\partial F(V)$, contradicting the fact that the restriction of
$F$ to $V-\{q\}$ is a local homeomorphism.

\vskip15pt

\section{A special homotopy and the proof of Theorem \ref{teorema 1}.}

Let $U\subset M$ be a neighborhood of $q$ on which a $C^1$ normal
(Gauss) map $\xi:U\to S^2$ is defined. Choose $a\in S^2$ such that
 $0<\langle a, \xi(q)\rangle <1$.
Shrinking $U$, one may assume that
$0<\langle a,
\xi(p)\rangle <1$ if $p\in U$.

Denote by $\alpha$ the second fundamental form of $M$, and
consider the height function $f:M\to \mathbb R$, $f(x)=\langle x,
a\rangle$. In particular, the (intrinsic) gradient and Hessian of
$f$ satisfy
\begin{eqnarray} \label{equation 1}
\nabla f(p)=a-\langle a, \xi(p)\rangle\xi(p),
\;\;\;\;\;\text{Hess} f(p)(v,v)=\langle \xi(p),a\rangle
\alpha(v,v), \;\; v\in T_pM.
\end{eqnarray}

The first equation expresses  the gradient of the restriction as the orthogonal projection of the space gradient  into the tangent space.  The  formula  for the Hessian of the restriction of a function  is  standard in submanifold geometry,  and can be found,  say, in (\cite {D}, p. 46).

Writing $H_f(p)$ and $A(p)$ for the linear endomorphisms
associated to the quadratic forms $\text{Hess} f(p)$ and
$\alpha(p)$ on $T_pM$,  respectively, it is clear from (\ref{equation 1}) that
the indices of the continuous symmetric  tensors $A$ and $H_f$
satisfy $j(H_f)=j(A)$.

Using the fact that $\det A $, being the Gaussian curvature, is negative  away from $q$ by hypothesis, it is easy to see
\begin{eqnarray*}
A(p)-\frac{t}{2}(\text{tr}A(p))I,\;\; 0\leq t \leq 1,
\end{eqnarray*}
is a homotopy through invertible maps  whenever  $p\neq q$.
Indeed, in a diagonalizing basis at $p$, the operator above has the matrix representation
$$\left[\begin{matrix}\lambda(p)(1-\frac{t}{2})-\mu(p) \frac{t}{2}&0\\0& \mu(p)(1-\frac{t}{2})-\lambda(p)\frac{t}{2}\end{matrix}\right]. $$
Since $t\in [0,1]$ and the principal curvatures satisfy  $\lambda(p)>0>\mu(p)$ if $p\neq q$, it is clear that the diagonal elements are non-zero.

It follows from the invariance of the degree under homotopies which  do not introduce further zeros that,
 for every continuous  non-vanishing  vector field $\eta$ on $U-\{q\}$,
$$j(A\eta)=j\big((A-\frac{1}{2}(\text{tr}A)I)\eta\big).$$ Hence, by
(\ref{equation 1}), Theorem  \ref{FormInd} and the fact that $\langle
\xi(p), a\rangle\neq 0$ for all $p\in U$,
\begin {eqnarray}\label{second index}
2j(A)=j(A\eta)+j(\eta)=j(H_f\eta)+j(\eta).
\end{eqnarray}
Applying (\ref{second index}) with $\eta=\nabla f$ (which is a
permissible choice, since $0<\langle a, \xi(p)\rangle <1 $ implies $\nabla
f(p)\neq 0$),  and using the general formula
\begin{eqnarray} \label{Hess}H_{\psi}\nabla \psi=
\nabla(\frac{1}{2}|\nabla \psi|^2), \end{eqnarray}
which is valid on any Riemannian manifold, with $\psi=f$,  one has
\begin{eqnarray} \label {third index}
2j(A)=j(H_f\nabla f)+j(\nabla f)=j\big(\nabla (\frac{1}{2}|\nabla
f|^2)\big)+j(\nabla f).
\end{eqnarray}
Manifestly, (\ref{third index}) provides a formula for the index
of the umbilic $q$ (i.e., the index of the tensor field $A$) in terms of the indices of two gradient
fields.

To conclude the proof of the theorem, we must show  that $j(A)$ is
non-positive. If $q$ is not an umbilical (planar) point, the
distributions of principal directions extend continuously across
$q$,  in which case $j(A)$ vanishes.  Hence, we may assume that
$q$ is an umbilic, that is, $\lambda(q)=\mu(q)=0$. In particular,
$A(q)=0$.

Since $\nabla f(p)\neq 0$ for $p\in U$, the term $j(\nabla f)$ in
(\ref{third index}) vanishes.  Hence, it remains to argue that
$j(\nabla (\frac{1}{2}|\nabla f|^2))\leq 0$. According to Lemma
\ref{IndGrad},   one must show that  the function
$h:=\frac{1}{2}|\nabla f|^2$ has an isolated critical point at
$q$, which is neither a local maximum nor a local minimum. From
(\ref{equation 1}) and (\ref {Hess}), one has $\nabla
h(p)=\nabla(\frac{1}{2}|\nabla f|^2)(p)=\langle a, \xi(p)\rangle
A(p)\nabla f(p)$. Since $A(q)=0$, the point  $q$ is critical for
$h$. In order to see that no point $p$ in $U-\{q\}$ is critical
for $h$, observe that $\langle a, \xi(p)\rangle >0$, $\nabla
f(p)\neq 0$ and $A(p)$ is invertible (since the Gaussian curvature
$\text{det}\; A(p)$ is negative for $p\neq q$, by hypothesis).
Hence $q$ is an isolated critical point of $h$.

We now proceed to show that $h$ has neither a local maximum nor a
minimum at $q$. A direct calculation gives $\sqrt {2h(p)}=|\nabla
f(p)|= \sin\theta (p)$, where $\theta(p)\in [0, \pi]$ is the angle
between the vectors  $a$ and  $\xi(p)$. Since
 $0<\langle a, \xi(q)\rangle <1$, it
follows that $\theta (q)\in (0, \frac{\pi}{2})$. Therefore,  in
order to show that $h$  does not have an extremum  at $q$, it suffices to argue that the Gauss
map $\xi:U\to S^2$ is open. But this is a consequence of Lemma
\ref{OpenMap} and the inverse function theorem, since the
Jacobian determinant of the normal  map is the Gaussian curvature which, by
hypothesis, is negative  away from  $q$.

\vskip20pt

\vskip10pt

$$
\begin{array}{lcccccccccl}
\text{Francisco Fontenele}            &&&&&&&&& & \text{Frederico Xavier}\\
\text{Departamento de Geometria}      &&&&&&&&& & \text{Department of Mathematics}\\
\text{Universidade Federal Fluminense}&&&&&&&&& & \text{University of Notre Dame}\\
\text{Niter\'oi, RJ, Brazil}          &&&&&&&&& & \text{Notre Dame, IN, USA}\\
\text{fontenele@mat.uff.br}           &&&&&&&&& & \text{fxavier@nd.edu}\\
\end{array}
$$

\end{document}